\documentclass{amsart}

\title[Representation as sum of primes]{Mean representation number of
integers as the sum of primes}

\author[G. Bhowmik]{Gautami Bhowmik}
\address{Universit\'e de Lille 1\\ 
Laboratoire Paul Painlev\'e UMR CNRS 8524\\
59655 Villeneuve d'Ascq Cedex\\
France}
\email{bhowmik@math.univ-lille1.fr}

\author[J.-C. Schlage-Puchta]{Jan-Christoph Schlage-Puchta}
\address{Universiteit Gent\\ 
Department of Pure Mathematics and Computer Algebra\\
Krijgslaan 281\\
Gebouw S22\\
9000 Gent\\
Belgium}

\email{jcsp@cage.ugent.be}

\newtheorem{theo}{Theorem}[section]

\newtheorem{lemm}{Lemma}

\begin{document}
\begin{abstract}
Assuming the Riemann Hypothesis we obtain asymptotic estimates for the 
mean value of the number of representations of an integer as a sum of
two primes. By proving a corresponding $\Omega$-term, we show that
our result is essentially the best possible. 
\end{abstract}

\maketitle

\section{Introduction and Results}
When studying the Goldbach conjecture
that every  even integer larger than $2$ is the sum of two primes
it is natural to consider the corresponding problem for the von
Mangoldt function $\Lambda$. Instead of showing that an even integer
$n$ is the sum of two primes, one aims at showing that   
$G(n)=\sum_{k_1+k_2=n}\Lambda(k_1)\Lambda(k_2)$ is sufficiently
large, more precisely, $G(n)>C\sqrt{n}$ implies the Goldbach
conjecture. It is known since long that this result is 
true for almost all $n$. It is easy to see that if $f$ is an
increasing function such that the Tchebychev function
$\Psi(x)=x+\mathcal{O}(f(x))$, then the 
mean value of $G(n)$ satisfies the relation
\[
\sum_{n\leq x}G(n)=x^2/2+O(xf(x)).
\]
If we consider the contribution of only one zero of the Riemann zeta function
$\zeta$, an error term of size $\mathcal{O}(f(x)^2)$ appears, which, under the 
current knowledge on zero free regions of $\zeta$, would not be significantly 
better than  $\mathcal{O}(xf(x))$. 
 Fujii\cite{Fu1} studied the error term of this mean
value under the Riemann Hypothesis (RH) and obtained 
$$\sum_{n\le x}G(n)=x^2/2+\mathcal{O}(x^{3/2})$$
which he later improved \cite{Fu2} to
\begin{equation}
\label{eq:F2}
\sum_{n\le x}G(n)=x^2/2+H(x)+(\mathcal{O}(x\log x)^{4/3})
\end{equation}
with $H(x)=-2\sum_{\rho}\frac{x^{1+\rho}}{\rho(1+\rho)}$,
where the summation  runs over all non-trivial zeros of $\zeta$.
In fact,the oscillatory term $H(x)$ is present even without assuming
RH, however, it is necessary for the error estimate above.

In this paper we prove that
\begin{theo}
\label{thm:ErrorG2}
Suppose that the RH is true. Then
we have 
\[
\sum_{n\leq x} G(n) = \frac{1}{2}x^{2} + H(x) +
\mathcal{O}(x\log^5 x),
\]
and
\[
\sum_{n\leq x} G(n) = \frac{1}{2}x^{2} + H(x) + \Omega(x\log\log x).
\]
\end{theo}
This confirms a conjecture of Egami and
Matsumoto\cite[Conj.~2.2]{EgMat}. Recently, Granville\cite{Gran} used
(\ref{eq:F2}) to obtain new characterisations of RH. 
The innovation of the present work is the idea to use the
distribution of primes in short intervals to estimate exponential sums
close to the point 0. Note that using the generalised Riemann
Hypothesis one could similarly find bounds for the exponential sums in
question in certain neighbourhoods of Farey fractions. Such a bound,
for example,  fixes
a gap in the proof of \cite[Theorem~1C]{Gran}. This approach can
further be used to study the meromorphic continuation of the generating
Dirichlet-series $\sum G(n)n^{-s}$, as introduced by Egami and
Matsumoto\cite{EgMat}, a topic we deal with elsewhere\cite{GNB}.

The $\log$-power in the error term can be improved, but reaching
$\mathcal{O}(n\log^3 n)$ would probably require some new idea.

We would like to thank the referee for suggesting the use of
Lemma~\ref{Lem:Gal} below, which lead to a substantial improvement.
\par
\section{Proofs.}
To prove the first part of our theorem, we compute the sum using the
circle method. We use the standard notation.

Fix a large real number $x$, set $e(\alpha)=e^{2\pi i \alpha}$ and let
\begin{eqnarray*}
S(\alpha) & = & \sum_{n\leq x}\Lambda(n) e(\alpha n),\\
T_y(\alpha) & = & \sum_{n\leq y} e(\alpha n),\\
T(\alpha) & = & T_x(\alpha),\\ 
R(\alpha) & = & S(\alpha)-T(\alpha).
\end{eqnarray*}
The following is due to Selberg\cite[eq. (13)]{Sel}.
\begin{lemm}
\label{Lem:Selberg}
Assuming RH we have 
\[
\int_1^x|\Psi(t+h)-\Psi(t)-h|^2\;dt \ll xh\log^2 x.
\]
\end{lemm}
The following result is due to Gallagher, confer\cite[Lemma~1.9]{Mon}
and put $T=y^{-1}, \delta=y/2$.
\begin{lemm}
\label{Lem:Gal}
Let $c_1, \ldots, c_N$ be complex numbers, and set $S(t)=\sum_{n=1}^N
c_n e(tn)$. Then
\[
\int\limits_{-1/y}^{1/y}|S(t)^2dt\ll
y^{-2}\int_{-\infty}^\infty|A(x)|^2dx.
\]
where
\[
A(x) = \underset{|n-x|\leq y/4}{\sum\limits_{n\leq N}}c_n.
\]
\end{lemm}
Our main technical result is the following.
\begin{lemm}
\label{Lem:L2local}
Suppose the RH. Then we have for $y\leq x$ the estimate
\[
\int\limits_{-y^{-1}}^{y^{-1}} |R(\alpha)|^2 d\alpha \ll
\frac{x}{y}\log^4 x.
\]
\end{lemm}
\begin{proof}
We put $N=x$ and $c_n=\Lambda(n)-1$ into
Lemma~\ref{Lem:GalSieve}. Putting
\[
B(t)=\underset{t<n\leq t+y/2}{\sum\limits_{n\leq x}}c_n
\]
we obtain
\[
\int\limits_{-y^{-1}}^{y^{-1}} |R(\alpha)|^2 d\alpha \ll
y^{-2}\int\limits_{-\infty}^\infty |B(t)|^2 dt =
y^{-2}\int\limits_{-y/2}^N |B(t)|^2 dt. 
\]
In the range $-y/2<t<0$ we have
\[
\int\limits_{-y/2}^0|B(t)|^2dt = \int\limits_0^{y/2}|\Psi(t)-[t]|^2dt
\ll y^2\log^4 y.
\]
For $0\leq t\leq x-y/2$ we have
$B(t)=\Psi(t+y/2)-\Psi(t)-y/2+\mathcal{O}(1)$, thus we can apply
Lemma~\ref{Lem:Selberg} to obtain
\[
\int\limits_0^{x-y/2} |B(t)|^2 dt \ll x +
\int\limits_0^{x-y/2}|\Psi(t+y/2)-\Psi(t)-y/2|^2dt \ll xy\log^2 x.
\]
Finally for $x-y/2\leq x\leq N$ we have
$B(x)=\Psi(x)-\Psi(t)-(x-t)+\mathcal{O}(1)$. The RH being equivalent
to $\Psi(x)=x+\mathcal{O}(x^{1/2}\log^2 x)$, this implies $B(x)\ll
x^{1/2}\log^2 x$, and therefore
\[
\int\limits_{x-y/2}^x |B(t)|^2 dt \ll xy\log^4 x.
\]
Collecting our estimates our claim follows.
\end{proof}

Note that no non-trivial unconditional version of
Lemma~\ref{Lem:L2local} can be proven without better understanding of
the zeros of the Riemann $\zeta$-function, since the existence of a
single zero close to 1 would already blow up the left-hand side.

Writing $S^2(\alpha)$ as $(T(\alpha)+R(\alpha))^2$ we have
\[
\sum_{n\leq x} G(n) = \int\limits_0^1
T(-\alpha)S^2(\alpha)d\alpha = \frac{1}{2}x^2 + 2\int\limits_0^1
|T(\alpha)|^2R(\alpha)d\alpha 
 + \int\limits_0^1
T(-\alpha)R^2(\alpha)d\alpha + \mathcal{O}(x).
\]
We claim that the second term yields $H(x)$, and the last one an error
of admissible size. In fact, the second term can be written
as
\[
2\int_0^1|T(\alpha)|^2S(\alpha)d\alpha-2\int_0^1|T(\alpha)|^2T(\alpha)d\alpha
 = 2\sum_{n\leq x}(\Lambda(n)-1)([x]-n) = 2\sum_{n\leq x-1}(\Psi(n)-n).
\]
We now insert the explicit formula for $\Psi(n)$, and replace the sum
over $n$ by an integral to find that the second term is indeed
$H(x)+\mathcal{O}(x)$. 

We now consider the third term.
We split the integral into an integral over
$[-x^{-1}, x^{-1}]$ and integrals of the form $[2^kx^{-1},
  2^{k+1}x^{-1}]$. On each interval we bound $T(\alpha)$ by
$\min\{x, \frac{1}{\|\alpha\|}\}$,
where $\|\alpha\|$ is the distance of $\alpha$ to the nearest integer,
and $R(\alpha)$ using  Lemma~\ref{Lem:L2local}. 
For the first interval this yields
\[
\int\limits_{-x^{-1}}^{x^{-1}} T(-\alpha)R^2(\alpha)d\alpha
\ll x\int\limits_{-x^{-1}}^{x^{-1}} R^2(\alpha)d\alpha\ll x\log^4 x,
\]
while for the other intervals we obtain
\begin{multline*}
\int\limits_{2^kx^{-1}}^{2^{k+1}x^{-1}} T(-\alpha)R^2(\alpha)d\alpha
\ll 2^{-k}x\int\limits_{2^kx^{-1}}^{2^{k+1}x^{-1}}
R^2(\alpha)d\alpha\\
\ll 2^{-k} x\frac{x}{2^{-k}x}\log^4 x
\ll x\log^4 x
\end{multline*}
There are $\mathcal{O}(\log x)$ 
summands, 
hence, the contribution of
$R^2$ to the whole integral is $\mathcal{O}(x\log^5 x)$, and the first
part of our theorem is proven. 

We now turn to the proof of the $\Omega$-result. To do so we show that
$G(n)=\Omega(n\log\log n)$, hence, the left hand side of (\ref{eq:F2}) has
jumps of order $\Omega(n\log\log n)$. Since $x^2/2$ and $H(x)$ are continuous,
the error term cannot be $o(x\log\log x)$. By considering the average
behaviour of $H(n)-H(n-1)$, one can even show that the error term is of order
$\Omega(x\log\log x)$ for integral $x$, however, we will only do the easier
case of real $x$ here.

The idea of the proof is that if an $n$ is divisible by many small primes,
then $G(n)$ should be large. Let $q_1$ be the exceptional modulus for which a
Siegel-zero for moduli up to $Q$ might exist, and $p_1$ be some prime divisor
of $q_1$. For the sake of determinacy we put $p_1=2$, if no Siegel zero exists.
We now use the following result due to Gallagher\cite[Theorem~7]{Gal}.

\begin{lemm}
\label{Lem:GalSieve}
We have
\[
\left|x-\sum_{x\leq n\leq x+h}\Lambda(n)\right|+\sum_{1<q\leq
  Q}{\sum_\chi}^*\left|\sum_{x\leq n\leq x+h}\Lambda(n)\chi(n)\right| \ll
h\exp\left(-c\frac{\log x}{\log Q}\right),
\]
provided that $x/Q\leq h\leq x$, $\exp(\log^{1/2} x)\leq Q\leq x^c$, $c$ is an
absolute positive constant, $\sum^*$ denotes summation over primitive
characters modulo $q$, and if there exists an exceptional character, for
which a Siegel zero exists, this character has to be left out of the summation.
\end{lemm}

We put $Q=q=\prod_{p<h, p\neq p_1} p$. Then all characters $\chi$ modulo $q$
is induced by some primitive character $\chi'$ modulo $q'\leq q$, and
\[
\left|\sum_{x\leq n\leq x+h}\Lambda(n)\chi(n) - \sum_{x\leq n\leq
  x+h}\Lambda(n)\chi'(n)\right| \leq \sum_{d|q} \Lambda(d)\leq \log q,
\]
which is negligible. Hence, it follows from Lemma~\ref{Lem:Gal}, that
\[
\left|x-\sum_{x\leq n\leq 2x}\Lambda(n)\chi_0(n)\right| +
\underset{\chi\neq\chi_0}{\sum_{\chi\pmod{q}}}\left|\sum_{x\leq n\leq 
  2x}\Lambda(n)\chi(n)\right| \leq \frac{x}{2},
\]
where $\chi_0$ is the principal character, provided that $q<x^{c'}$ for some
absolute constant $c'$. It follows that for $(a, q)=1$ we have
\[
S(x, q, a) := \underset{n\equiv a\pmod{q}}{\sum_{n\leq x}}\Lambda(n)\geq\frac{x}{\varphi(q)}.
\]
Now
\[
\underset{q|n}{\sum_{n\leq 4x}} G(n)\geq \sum_{(a, q)=1}S(x, q, a)S(x, q, q-a) \geq\frac{x^2}{4\varphi(q)}.
\]
On the left we take the average over $\ll\frac{x}{q}$ integers, hence, we obtain
\[
\max_{n\leq 4x} G(n) \gg \frac{x}{2\varphi(q)}= (1-p_1^{-1})\prod_{p\leq
h}(1-p^{-1})^{-1} x \gg x\log\log x,
\]
and our claim follows.
\par

\end{document}